\documentstyle[12pt]{article}
 \begin{document}
 \baselineskip=15pt
\newcommand{\be}{\begin{equation}}
\newcommand{\ee}{\end{equation}}
 \newcommand{\ekq}{e^{2\pi i\frac{k}{q}}}
 \newcommand{\rf}{Ramanujan - Fourier~}
 \newcommand{\ekqn}{e^{2\pi i\frac{k}{q}n}}
 \newcommand{\ekqm}{e^{-2\pi i\frac{k}{q}}}
 \newcommand{\ekqmn}{e^{-2\pi i\frac{k}{q}n}}
 \newcommand{\wk}{Wiener - Khintchine formula~}
 \newcommand{\ekqp}{e^{2\pi i\frac{k'}{q'}}}
 \newcommand{\di}{\displaystyle}
 \begin{center}
{\Large {\bf Ramanujan - Fourier Series and \\
the Density of Sophie Germain Primes}}\\
 \vspace{1cm}
H. Gopalkrishna Gadiyar and R. Padma \\
AU-KBC Research Centre, M. I. T. Campus of Anna University, Chromepet, Chennai 600 044, India\\
E-mail: \{gadiyar, padma\}@au-kbc.org
\end{center}

\noindent {\bf Abstract.} A prime $p$ is called Sophie Germain prime if $2p+1$ is also prime. A formula for the density of such primes is given in a more general setting using a new approach. This method uses the \rf series for a modified van Mangoldt function. The proof remains heuristic as interchange of certain limits has not been justified. Experimental evidence using computer calculations is provided for the plausibility of the result.
\\
\\

\noindent {\bf Introduction.} Ramanujan's work has been a source of inspiration to researchers in the domain of experimental mathematics [Borwein, Borwein, Bailey 1989]. Ramanujan himself used extensive numerical calculations before making conjectures and one wonders what he would have done with the present technology of computers. Fourier series are an extensively used tool across disciplines as a means of extracting useful information from seemingly disorganized data. In this paper we wish to bring to the attention of the experimental mathematics community the Ramanujan-Fourier series which are useful for analyzing arithmetical functions. We use this tool to give a heuristic ``proof" (which can be made rigorous if certain limits can be interchanged) and give numerical evidence for the formula derived.

A positive integer $p$ is called a Sophie Germain prime if both $p$ and $2p+1$ are primes, (2,5), (3,7), (5,11), (11,23) for example. The question is: Are there infinitely many Sophie Germain primes? The conjecture is: Yes. The density of Sophie Germain primes has been conjecured in [Hardy, Littlewood 1922] following probabilistic reasoning of Brun. Sophie Germain primes are of great interest in recent times after the famous AKS algorithm for primality testing [Agrawal, Kayal, Saxena, 2002]. If the conjecture about the density of Sophie Germain primes is true, then the complexity of the AKS algorithm can be brought down to $O(\log ^6n)$. Sophie Germain primes are the most sought after primes for the RSA algorithm as they are robust against the Pollard's $p-1$ method of factoring [Stinson 1995]. 

In [Gadiyar, Padma 1999], we gave a heuristic proof of the twin prime conjecture using the \rf series. In this paper, we show that the ``proof" could be extended to a more general class of primes in which Sophie Germain prime is a special case. Unlike the earlier probabilistic arguments, our approach is analytic. We show that the conjecture regarding the density of prime pairs $(p, p')$ satisfying the condition $ap' - bp~=~l$, where $a$ and $b$ are positive integers, is true upto interchange of certain limits. 

We use the following two tools. 1. The \rf series for ${\di \Lambda_1(n)= \frac{\phi(n)}{n}} \Lambda (n)$ where $\Lambda (n)$ is the von Mangoldt function defined as $\log p$ if $n=p^k$, where $p$ is prime and $k$ any positive integer and is equal to $0$ otherwise. $\phi (n)$ denotes the number of integers less than or equal to n and relatively prime to $n$. 2. Carmichael's formula for getting the Ramanujan - Fourier coefficients for arithmetical functions.  

Numerical evidence is given for various choices of $a$ and $b$ which shows remarkable accuracy of the conjecture.

The paper consists of the following sections. In  Section 1, we state the main ``result", in Section 2, we describe the \rf series and some properties of the Ramanujan sum, in Section 3, we give the ``proof" of main result and in Section 4, we give the numerical evidence of our main result for various choices of $a$ and $b$.

\noindent {\bf 1. Main Conjecture.} Let $a$, $b$ and $l$ be positive integers, where $(a,b)=1$.  Let us ask the question whether there are infinitely many prime pairs $(p,p')$ satisfying the equation $ap' - bp~=~l$? 
In [Hardy, Littlewood, 1922], Hardy and Littlewood conjectured the following. Let $\pi_{(a,b,l)}(N)$ denote the number of prime pairs $(p,p')$ satisfying the condition $ap' - bp~=~l$ such that $p<N$. Then
\begin{equation}
\pi_{(a,b,l)}(N) = o\left (\frac{N}{\log ^2N}\right )
\end{equation}
unless $(l,a)=1$, $(l,b)=1$, and just one of $a,b,l$ is even. But if these conditions are satisfied then
\begin{equation}
\pi_{(a,b,l)}(N) \sim \frac{2C}{a}\frac{N}{\log ^2N} \prod_{\stackrel{p | abl}{p >2}} \left ( \frac{p-1}{p-2}\right ) \, ,
\end{equation}
where 
\begin{equation}
C=\prod_{p > 2} \left ( 1-\frac{1}{(p-1)^2} \right ) \, ,
\end{equation}
and $p$ denotes a prime. 

\noindent {\bf Main Result.} Let $\Psi_{(a,b,l)}(N)=\sum_{n \le N} \Lambda_1(n)\Lambda_1\left (\frac{bn+l}{a}\right )$. Then, upto interchange of certain limits,
\be
\lim_{N \rightarrow \infty} \frac{1}{N}\Psi_{(a,b,l)}(N)  = \left \{ \begin{array}{ll}\frac{2C}{a} \prod_{\stackrel{p | abl}{p >2}} \left ( \frac{p-1}{p-2}\right ) \,,& \rm{if~}(a,l)=(b,l)=1 \\
~\rm{and~exactly~one}&\rm{of~}a,b,l \rm{~is~even,}\\
~&~\\
0, &\rm{otherwise.} \end{array} \right . 
\ee
Note that (1) and (2) follow immediately from (4), see for example [Gadiyar, Padma 1999].

The twin prime problem corresponds to the case $a=b=1$ and $l=2$. The Sophie Germain prime problem corresponds to the case $a=1$, $b=2$ and $l=1$. In [Gadiyar, Padma 1999], we showed that the twin prime problem is related to autocorrelation and hence to the \wk which is used in probability and electrical engineering. We show in Section 3 that the heuristic proof that was given in that paper can be extended to prove (4). To do this, we need a short description of \rf series which we give in the next section.

\noindent {\bf 2. Ramanujan - Fourier Series.} Ramanujan in [Ramanujan 1918] showed that many important arithmetical functions ($a(n)$) have an expansion of the form
 \begin{equation}
 a(n) ~=~ \sum_{q=1}^{\infty} a_q c_q(n) \, ,
 \end{equation}
 where
 \begin{equation}
 c_q(n) ~=~ \sum_{\stackrel{k=1}{(k,q)=1}}^q \ekqn \, ,
\end{equation}
is called the Ramanujan sum and the $a_q$'s are known as the Ramanujan-Fourier coefficients. He obtained such expansions for $d(n)$, $\sigma (n)$, $\phi(n)$ and so on where $d(n)$ denotes the number of divisors of $n$ and $\sigma (n)$ denotes the sum of divisors of $n$. In [Hardy 1921], Hardy proved that the Ramanujan sum is a multiplicative function of $q$, that is,
\begin{equation}
c_{q q'}=c_q(n) c_{q'}(n) \, \rm{~if~} (q,q')=1 \, ,
\end{equation}
using which he obtained the \rf expansion of $\Lambda _1(n)$.
\begin{equation}
\Lambda _1 (n) = \sum_{q=1}^\infty \frac{\mu(q)}{\phi(q)}
 c_q(n) \, ,
\end{equation}
where $\mu (q)$ is the M\"{o}bius function defined as follows:
\begin{equation}
 \mu (q)~=~ \left \{ \begin{array}{ll} (-1)^k \, , & {\rm if~}q~=~
 p_1p_2...p_k,~ p_i{\rm 's ~are~ distinct~ primes \, ,}\\
 0 \, , & {\rm otherwise.}
 \end{array} \right . 
\end{equation}
If $p$ is prime, then
\begin{equation}
c_p(n)~=~ \left \{ \begin{array}{ll} -1 \, , &{\rm ~if ~} p \mid\!\!\!\!/ n \, ,
 \\
 p-1, & {\rm ~if~} p | n \, ,
 \end{array} \right .
\end{equation}
where $a | b$ means $a$ divides $b$ and $a \mid\!\!\!\!/ b$ means $a$ does not divide $b$.

Neither Ramanujan nor Hardy gave a formula for finding the \rf coefficients which are the back bone of Fourier analysis. This was done later by Carmichael [Carmichael 1932]. Let $M(f)$ denote the mean value of an arithmetical function $f$, that is,
 \be
 M(f) ~=~ \lim_{N \rightarrow \infty} \frac{1}{N} \sum_{n \le N} f(n) \, .
 \ee
 For $1 \le k \le q, ~ (k,q)=1,$ let
 $e_{\frac{k}{q}}(n)=e^{2\pi i \frac{k}{q} n}$, ($n \in \cal N$). If
 $a(n)$ is an arithmetical function with expansion (5)
 then
 \be
 a_q ~=~ \frac{1}{\phi (q)} M(a ~c_q) ~= \frac{1}{\phi (q)} ~\lim_{N \rightarrow
 \infty} \frac{1}{N}
 \sum_{n \le N} a(n) c_q(n)   ~.
 \ee
 Also,
 \be
 M(e_{\frac{k}{q}}  ~\overline{e_{\frac{k'}{q'}}}) = \left \{ \begin{array}{ll}
 1 \, , &{\rm ~if~} \frac{k}{q}~=~ \frac{k'}{q'} \, ,\\
 0 \, , & {\rm ~if~} \frac{k}{q} ~\ne ~\frac{k'}{q'} \, . 
 \end{array} \right .
 \ee

\noindent {\bf 3. ``Proof" of Main Result.} Using the results in Section 2, we give the heuristic proof of (4). Note that $\Lambda(\frac{bn+l}{a})$ is defined if and only if $a | (bn+l)$. Using the identity
\be
\frac{1}{a} \sum_{j=0}^{a-1} e^{2\pi i\frac{j}{a}m}=\left \{ \begin{array}{ll}
 1 \, , &{\rm ~if~}a | m \, ,\\
 0 \, , & {\rm ~if~}a \mid\!\!\!\!/ m \, , 
 \end{array} \right .
\ee
we write 
\begin{eqnarray}
\lim_{N \rightarrow \infty} \frac{1}{N}\Psi_{(a,b,l)}(N)&=&\lim_{N \rightarrow \infty}\frac{1}{N} \sum_{n \le N} \Lambda_1(n)\Lambda_1\left (\frac{bn+l}{a}\right ) \frac{1}{a} \sum_{j=0}^{a-1} e^{-2\pi i\frac{j}{a}(bn+l)} \nonumber \\
~&\stackrel{?}{=}&\left (\frac{1}{a} \sum_{j=0}^{a-1}\sum_{q=1}^{\infty}\sum_{\stackrel{k=1}{(k,q)=1}}^q \sum_{q'=1}^{\infty}\sum_{\stackrel{k'=1}{(k',q')=1}}^{q'}\frac{\mu(q)}{\phi(q)}\frac{\mu(q')}{\phi(q')}e^{-2\pi i(\frac{k'}{q'}+j)\frac{l}{a}} \right )\nonumber\\
~&&~~~~~\left (\lim_{N \rightarrow \infty}\frac{1}{N} \sum_{n \le N}e^{2\pi i \left (\frac{k}{q}- \frac{k'}{q'} \frac{b}{a}-j \frac{b}{a} \right )n}\right ) \\
~&=&\frac{1}{a} \sum_{j=0}^{a-1}\sum_{q=1}^{\infty}\sum_{\stackrel{k=1}{(k,q)=1}}^q \sum_{q'=1}^{\infty}\sum_{\stackrel{k'=1}{\stackrel{(k',q')=1}{\frac{k}{q}=\left ( \frac{k'}{q'}+j \right ) \frac{b}{a}}}}^{q'}\frac{\mu(q)}{\phi(q)}\frac{\mu(q')}{\phi(q')}e^{-2\pi i(\frac{k'}{q'}+j)\frac{l}{a}} \nonumber \\
~&&~\\
~&=& S {\rm ~(say),} \nonumber
\end{eqnarray}
where we have used (8), freely interchanged the sums and limits to obtain (15), and then used (13) to get (16). We prove that $S$ is equal to the R. H. S. of (4) as a lemma.

\noindent {\bf Lemma.} If $(a,b)=1$, then
\be
S=\left \{ \begin{array}{ll}\frac{2C}{a} \prod_{\stackrel{p | abl}{p >2}} \left ( \frac{p-1}{p-2}\right ) \,,& \rm{if~}(a,l)=(b,l)=1, \\
~\rm{and~exactly~one}&\rm{of~}a,b,l \rm{~is~even}\\
~&~\\
0, &\rm{otherwise.} \end{array} \right . 
\ee
{\bf Proof.} Since $(k',q')=1$ and $0\le j \le a-1$, 
\begin{eqnarray}
\frac{k}{q}&=&\left ( \frac{k'}{q'}+j \right ) \frac{b}{a} \nonumber \\
~&=&\frac{k'+jq'}{q'a}~b \nonumber \\
&=&\frac{k_1}{q'a}~b, {\rm ~where~}1\le k_1 \le q'a {\rm ~and~}(k_1,q')=1.
\end{eqnarray}
Now, (18) can happen if and only if 
\be
qd_2=q'd_1 \, ,
\ee
for some divisors $d_1$ and $d_2$ of $a$ and $b$ respectively and  
\begin{eqnarray}
k&=&\frac{k_1}{q'a}bq \nonumber \\
~&=&\frac{k_1\frac{b}{d_2}qd_2}{q'd_1 \frac{a}{d_1}}  \nonumber \\
~&=&k_1 \frac{ \frac{b}{d_2}}{\frac{a}{d_1}} \, ,
\end{eqnarray}
is an integer. Since $(a,b)=1$, this can happen if and only if $\frac{a}{d_1}$ divides $k_1$. Also, from (19), $d_1 |qd_2$ and since $(a,b)=1$, $d_1|q$. So we write $q=d_1q_1$ where $q_1\ge 1$ is an integer. Similarly $q'=d_2q_2$ where $q_2\ge 1$ is an integer. Thus from (19), $q_1=q_2$.
Also $\mu(q)=\mu(d_1q_1) \neq 0$ if and only if $(d_1,q_1)=1$. Similarly we have $(d_2,q_2)=1$. 

Let us write $k_1= \frac{a}{d_1} k_2$. Since $(k,q)=1$ and $d_1|q$, we have $(d_1,k)=1$ and from (20), $k=k_2 \frac{b}{d_2}$ and thus $(k_2,d_1)=1$. So we can write
\be
S=\frac{1}{a} \sum_{d_1|a} \sum_{d_2|b} \sum_{\stackrel{q_2=1}{\stackrel{(d_1,q_2)=1}{(d_2,q_2)=1}}}^{\infty} \frac{\mu(d_1)}{\phi(d_1)}\frac{\mu(d_2)}{\phi(d_2)}\frac{\mu^2(q_2)}{\phi^2(q_2)}\sum_{\stackrel{k_2=1}{(k_2,q_2d_1d_2)=1}}^{q_2d_1d_2} e^{-2\pi i \frac{k_2}{q_2d_1d_2}l}\, .
\ee
That is,
\be
S=\frac{1}{a} \sum_{d_1|a} \sum_{d_2|b} \sum_{\stackrel{q_2=1}{\stackrel{(d_1,q_2)=1}{(d_2,q_2)=1}}}^{\infty} \frac{\mu(d_1)}{\phi(d_1)}\frac{\mu(d_2)}{\phi(d_2)}\frac{\mu^2(q_2)}{\phi^2(q_2)}c_{(q_2d_1d_2)}(l) \, ,
\ee
by the definition of the Ramanujan sum (6).
Now by the multiplicative property (7) of the Ramanujan sum, since $q_2,d_1,d_2$ are pairwise relatively prime, 
\begin{eqnarray}
S&=&\frac{1}{a} \sum_{d_1|a} \sum_{d_2|b} \sum_{\stackrel{q_2=1}{\stackrel{(d_1,q_2)=1}{(d_2,q_2)=1}}}^{\infty} \frac{\mu(d_1)}{\phi(d_1)}\frac{\mu(d_2)}{\phi(d_2)}\frac{\mu^2(q_2)}{\phi^2(q_2)} c_{q_2}(l)c_{d_1}(l)c_{d_2}(l) \nonumber \\
~&=&\frac{1}{a}\sum_{\stackrel{q_2=1}{\stackrel{(a,q_2)=1}{(b,q_2)=1}}}^{\infty}\frac{\mu^2(q_2)}{\phi^2(q_2)} c_{q_2}(l)\sum_{d_1|a}\frac{\mu(d_1)}{\phi(d_1)}c_{d_1}(l)\sum_{d_2|b}\frac{\mu(d_2)}{\phi(d_2)}c_{d_2}(l) \nonumber \\
~&=&\frac{1}{a}\sum_{\stackrel{q_2=1}{(ab,q_2)=1}}^{\infty}\frac{\mu^2(q_2)}{\phi^2(q_2)} c_{q_2}(l)\sum_{d_1|a}\frac{\mu(d_1)}{\phi(d_1)}c_{d_1}(l)\sum_{d_2|b}\frac{\mu(d_2)}{\phi(d_2)}c_{d_2}(l) \, ,
\end{eqnarray}
as $(a,b)=1$. Writing the series and sums in (23) as Euler products, we get
\be
S=\frac{1}{a}\prod_{p\mid\!\!\!/ab} \left ( 1+\frac{c_p(l)}{(p-1)^2} \right ) \prod_{p|a} \left(1-\frac{c_p(l)}{p-1} \right ) \prod_{p|b} \left(1-\frac{c_p(l)}{p-1} \right ) \, .
\ee
By the property (10) of the Ramanujan sum, if $u>1$ is an integer, then
\be
\prod_{p|u} \left(1-\frac{c_p(l)}{p-1} \right ) = \left \{\begin{array}{ll}
 \prod_{p|u} \left(\frac{p}{p-1} \right ) \, , &{\rm ~if~}p \mid\!\!\!\!/ l \, ,\\
 0 \, , & {\rm ~if~}p | l \, . 
 \end{array} \right .
\ee
Hence we will assume that $(a,l)=1$ and $(b,l)=1$ so that
\be
S=\frac{1}{a}\prod_{p\mid\!\!\!/abl} \left ( 1-\frac{1}{(p-1)^2} \right )\prod_{\stackrel{p\mid\!\!\!/ab}{p|l}} \left (\frac{p}{p-1}\right )\prod_{\stackrel{p|a}{p\mid\!\!\!/l
}} \left (\frac{p}{p-1}\right )\prod_{\stackrel{p|b}{p\mid\!\!\!/l
}} \left (\frac{p}{p-1}\right ) \, .
\ee
If none of $a,b,$ or $l$ is even, then the product
\be
\prod_{p\mid\!\!\!/abl} \left ( 1-\frac{1}{(p-1)^2} \right )=0 \, .
\ee
So we will assume that one of $a,b$ or $l$ is even. But $(a,b)=1$, $(a,l)=1$ and $(b,l)=1$ and therefore exactly one of $a,b$ or $l$ is even. We can therefore write the infinite product
\be
\prod_{p\mid\!\!\!/abl} \left ( 1-\frac{1}{(p-1)^2} \right )= \prod_{p>2} \left ( 1-\frac{1}{(p-1)^2} \right )\prod_{\stackrel{p|abl}{p>2}} \left (\frac{(p-1)^2}{p(p-2)} \right ) \, .
\ee
Thus $S=0$ unless $(a,l)=1$, $(b,l)=1$ and exactly one of $a,b$ or $l$ is even, but if these conditions are satisfied, then the value of S as given in (17) is got by simplifying (26) using (28).

\noindent{\bf 4. Experimental Evidence.} We give now the compelling numerical evidence of the main result (4) by varying $a$ and $b$. We have taken the value of $C \sim 0.660161816$ and Ratio is defined by (R.H.S. of (4))/${\di \frac{\Psi_{(a,b,l)}(N)}{N}}$.

\noindent{{\bf Example 1.}} We take $a=1,~b=2,~l=1$ which corresponds to Sophie Germain primes. In this case, the R. H. S. of (4) = $2C=1.320323632$.

\begin{center}
 {\bf Table 1}\\
 ~\\
 \begin{tabular}{|c|c|c|c|}
 \hline
 \hline
~&~&~&~\\
 N& $\Psi_{(1,2,1)}(N)$ & ${\di \frac{\Psi_{(1,2,1)}(N)}{N}} $ & Ratio\\
~&~&~&~\\
 \hline
 50000 & 66130.966133  & 1.322619  & 0.998264 \\
 100000 & 132886.401744 & 1.328864 & 0.993573 \\
 150000 & 200755.416380  & 1.338369  & 0.986517 \\
 200000 & 265612.706085  & 1.328064 & 0.994172  \\
 250000 & 331585.551940  & 1.326342  &0.995462 \\
 300000 & 394316.641234  & 1.314389  & 1.004515  \\
 350000 & 459668.599011 & 1.313339  & 1.00531 \\
 400000 & 521496.993567 & 1.303742 & 1.012718 \\
 450000 & 588393.432192 &1.307541  & 1.009776 \\
 500000 & 652614.182933  & 1.305228 & 1.011565 \\
 \hline
 \end{tabular}

\end{center}
\newpage
\noindent{{\bf Example 2.}} We take $a=1,~b=10,~l=1$. In this case, the R. H. S. of (4) = $\frac{8C}{3}=1.760431509$.

\begin{center}
 {\bf Table 2}\\
~ \\
 \begin{tabular}{|c|c|c|c|}
 \hline
 \hline
~&~&~&~\\
 N& $\Psi_{(1,10,1)}(N)$ & ${\di \frac{\Psi_{(1,10,1)}(N)}{N}} $ & Ratio\\
~&~&~&~\\
 \hline
 10000 & 17107.791529 & 1.710779  & 1.029023 \\
 20000 & 34210.057148 & 1.710503 & 1.029189 \\
 30000 & 51939.100560   & 1.731303  & 1.016824 \\
 40000 & 70219.348038   & 1.755484  & 1.002818   \\
 50000 & 89934.594398   & 1.798692   &0.978729  \\
 60000 & 106902.836342 & 1.781714  &0.988055  \\
 70000 & 123796.944818 &1.768528    & 0.995422  \\
 80000 & 141470.265879  & 1.768378 & 0.995506  \\
 90000 & 159287.348829 &1.769859   & 0.994673 \\
 100000 & 177824.093558  & 1.778241 & 0.989985 \\
 \hline
 \end{tabular}

\end{center}

\noindent{{\bf Example 3.}} We take $a=3,~b=5,~l=2$. In this case, the R. H. S. of (4) = $\frac{16C}{9}=1.173621006$.

\begin{center}
 {\bf Table 3}\\
 ~\\
 \begin{tabular}{|c|c|c|c|}
 \hline
 \hline
~&~&~&~\\
 N& $\Psi_{(3,5,2)}(N)$ & ${\di \frac{\Psi_{(3,5,2)}(N)}{N}} $ & Ratio\\
~&~&~&~\\
 \hline

 60000 & 69649.061665& 1.160837 & 1.011013\\
 120000 & 140371.214304 &1.169770 & 1.003292  \\
 180000 & 211924.646933& 1.177366  &0.996819 \\
 240000 & 282504.323361  &1.177106   & 0.997039 \\
 300000 & 355072.360724    & 1.183578&0.991587  \\
 360000 & 423152.712312 & 1.175427 & 0.998463 \\
 420000 & 496296.973007& 1.181662 &  0.993195  \\
 480000 & 568659.361599  &1.184709 & 0.990640  \\
 540000 & 642488.622118 &1.189796 &0.986405  \\
 600000 & 712048.221861 & 1.186749 & 0.988938 \\
 \hline
 \end{tabular}

\end{center}

\noindent {\bf 5. Conclusion.} If the step (15) could be proved rigorously which involves justification of interchange of certain limits, then a whole class of outstanding problems including the twin prime problem and the Sophie Germain prime problem could be solved completely. We may say that in a precise sense  the \rf series for the (refined) von Mangoldt function traps the fluctuations in the distribution of primes. It is hoped that the theory of \rf series could be developed to study various properties of arithmetical functions. Numerical agreement between conjecture and experiment means that this technique could become a common tool and lead to further developments in number theory.  
\\
\\
\noindent {\bf References.}

\begin{description}
\item{[Agrawal, Kayal, Saxena, 2002]} M. Agrawal, S. Kayal, N. Saxena, ``Primes is in P", Annals of Mathematics, 160 (2004), 781-793.
\item{[Borwein, Borwein, Bailey 1989]} J.M. Borwein, P.B. Borwein, and D. A. Bailey, ``Ramanujan, modular equations and pi or how to compute a billion digits of pi", MAA Monthly, 96 (1989), 201-219.
\item {[Carmichael 1932]} R. D. Carmichael,  ``Expansions of arithmetical functions in infinite series", Proc. London Math. Soc. (2) 34 (1932), 1-26.
\item{[Gadiyar, Padma 1999]} H. Gopalkrishna Gadiyar and R. Padma, ``Ramanujan-Fourier series, the Wiener-Khintchine formula and the distribution of prime pairs", Physica A 269(1999), 503-510.
\item {[Hardy 1921]} G. H. Hardy, ``Note on Ramanujan's trigonometrical function $c_q(n)$ and certain series of arithmetical functions",  Proc. Camb. Phil. Soc. 20 (1921), 263-271.
\item {[Hardy, Littlewood 1922]}  G. H. Hardy and J. E. Littlewood, ```Some problems of Partition Numerorum'; III: On the expression of a number as a sum of primes",  Acta Math. 44 (1922), 1-70.
\item {[Ramanujan 1918]} S. Ramanujan, ``On certain trigonometrical sums and their applications in the theory of numbers", Trans. Camb. Phil. Soc.
 22 (1918), 259-276.
\item{[Stinson 1995]} D. R. Stinson, ``Cryptography: Theory and Practice", CRC Press, Boca Raton, 1995.
\end{description}

\end{document}